\documentclass[12pt, twoside, leqno]{article}

\usepackage{amsmath,amsthm}

\usepackage{amssymb,latexsym}

\usepackage{enumerate}

\usepackage{supertabular}

\newtheorem{theorem}{Theorem}
\newtheorem{proposition}[theorem]{Proposition}

\numberwithin{equation}{section}

\def\Q{{\mathbb Q}}
\def\Z{{\mathbb Z}}

\def\C{{\mathbb{C}}}

\def\gcd{\mathrm{gcd}}
\def\al{{\alpha}}
\def\be{{\beta}}
\def\ga{{\gamma}}
\def\de{{\delta}}
\def\la{{\lambda}}
\def\th{{\theta}}
\def\ep{{\epsilon}}

\def\ida{{\mathfrak a}}
\def\idb{{\mathfrak b}}
\def\idp{{\mathfrak p}}
\def\ordp{\mathrm{ord}_p}

\def\Arg{{\mathrm{Arg}}}
\def\Log{{\mathrm{Log}}}
\def\bi{{\mathbf{i}}}

\def\proof{{\bf Proof}.\:}
\newcommand{\proofend}{\hspace*{1mm} \hfill{$\Box$}}

\newcommand{\Nm}[2]{\mathrm{N}_{#1/#2}}
\newcommand{\ideal}[1]{\langle #1\rangle}

\newcommand{\cnj}[2]{#1^{(#2)}}

\begin{document}

\baselineskip=17pt

\title{On the Diophantine Equation $x^{2}+5^{a}\cdot 11^{b}=y^{n} $}

\author{I.N.~Cang\"{u}l  
\and
M.~Demirci
\and
G.~Soydan 
\and
N.~Tzanakis 
}

\date{} 

\maketitle

\renewcommand{\thefootnote}{}



\renewcommand{\thefootnote}{\arabic{footnote}}

\setcounter{footnote}{0}

\begin{abstract}
We give the complete solution $(n,a,b,x,y)$ of the title equation when $\gcd(x,y)=1$, except for the case
when $xab$ is odd. Our main result is Theorem \ref{main result}.
\end{abstract}

\section{Introduction}

The literature on the exponential Diophantine equation 
\begin{equation} \label{eq general}
x^{2}+C=y^{n}\,,\quad x\geq 1, y\geq 1, n\geq 3 
\end{equation}
goes back to 1850 when Lebesque \cite{Lebesque} proved that the
equation (\ref{eq general}) has no solutions when $C=1$. The title equation is actually a
special case of the Diophantine equation $ay^{2}+by+c=dx^{n}$,  where $a,b,c$
and $d$ are integers, $a\neq 0$, $b^{2}-4ac\neq 0$,  $d\neq 0$,  which has
only a finite number of solution in integers $x$ and $y$ when $n\geq 3$;  see 
\cite{Landau}.  J.H.E.~Cohn \cite{Cohn1}, solved (\ref{eq general}) for most values of $C$
in the range $1\leq C\leq 100$. The equations $x^{2}+74=y^{5}$ and $x^{2}+86=y^{5}$ that are not
solved in that paper, were later solved by Mignotte and de Weger in \cite{Mignotte}, and the
remaining unsolved cases in Cohn's paper were solved by Bugeaud, Mignotte and Siksek in \cite{Bugeaud2}.

Upper bounds for the exponent $n$ can be obtained as an application of the work of 
B\'{e}rczes, Brindza and Hajdu \cite{Berczes2} and of Gy\H{o}ry \cite{Gyory}. These results are based
on the Theory of Linear Forms in Logarithms and the obtained upper bounds, though effective, are not explicit.

Recently, the case in which $C$ is a power of a fixed prime gained the interest of several authors. 
In \cite{Arif4}, Arif and Muriefah solve $x^{2}+2^{k}=y^{n}$ under certain assumptions. 
In \cite{Le1}, Le verifies a conjecture of J.H.E~Cohn saying that $x^2+2^k=y^n$ has no solutions with
even $k>2$ and $x$ odd, which was proposed in \cite{Cohn2}.
The equation $x^{2}+3^{m}=y^{n}$ is completely solved by Arif and Muriefah in \cite{Arif1} when $m$ is
odd and by Luca in \cite{Luca3} when $m$ is even.
Liqun solves the same equation independently in \cite{Liqun1} for both odd and even $m$.
All solutions of $x^{2}+5^{m}=y^{n}$ with $m$ odd are given by Arif and Muriefah in \cite{Arif2} and
with $m$ even by Muriefah in \cite{FSAbu7}.
Again, the same equation is independently solved by Liqun in \cite{Liqun2}.
In \cite{Arif3}, Arif and Muriefah give the complete solution of $x^{2}+q^{2k+1}=y^{n}$ for $q$ odd prime, 
$q\not\equiv 7 \pmod{8}$ and $n\geq 5$ prime to $6h$, where $h$ is the class-number of the number field
$\Q(\sqrt{-q})$. 
Luca and Togbe solve $x^{2}+7^{2k}=y^{n}$ in \cite{Luca4} and B\'{e}rczes and Pink \cite{Berczes} 
solve (\ref{eq general}) with $C=p^{2k}$, where $2\leq p<100$ is prime, $(x,y)=1$ and $n\geq 3$.

More complicated cases, in which $C$ is a product of at least two prime powers are considered in some 
recent papers. For example, the complete solution $(n,x,y)$ with $n\geq 3$ and $\gcd(x,y)=1$
of the equation (\ref{eq general}), when $C$ is one of $2^a 3^b$, $5^a 13^b$, $2^a 5^b 13^c$, 
$2^a 11^b$, $2^a 3^b 11^c$ is respectively given in \cite{Luca1}, \cite{FSAbu3}, \cite{Luca5}, \cite{CDLPS}, 
\cite{CDIPS}.
In \cite{Luca2} the equation (\ref{eq general}) with $C=2^a 5^b$ is solved when $n\in\{3,4,5,6,8\}$
and $\gcd(x,y)=1$.
In \cite{Pink} all the {\em non-exceptional solutions} (in the terminology of that paper) of 
the equation (\ref{eq general}) with $C=2^a 3^b 5^c 7^d$ are given (with $n\geq 3$). Note that finding all 
the {\em exceptional solutions} of this equation seems to be a very difficult task.
\\
A survey of many relevant results can be found in \cite{FSAbu1}.

In this paper, we study the equation
\begin{equation} \label{main eqn}
x^{2}+5^{a}\cdot 11^{b}=y^{n}\,,\quad x\geq 1,y\geq 1, (x,y)=1,\, n\geq 3,\, a\geq 0,\, b\geq 0\,.
\end{equation}
Our main result is the following.
\begin{theorem} \label{main result}
When $n=3$, the only solutions to the equation (\ref{main eqn}) are 
\begin{align}
 (a,b,x,y) = & (0,1,4,3),\,(0,1,58,15),\,(0,2,2,5),\,(0,3,9324,443),\,(1,1,3,4), \nonumber \\
             & (1,1,419,56),\,(2,3,968,99),\,(3,1,37,14),\,(5,5,36599,1226)\,, \label{complete solution}
\end{align}
and, consequently, $(a,b,x,y)=(1,1,3,2)$ is the only solution when $n=6$. \\
When $n=4$, the equation (\ref{main eqn}) has no solutions.\\
When $n\geq 5, n\neq 6$, the equation (\ref{main eqn}) has no solutions $(a,b,x,y)$ with $ab$ odd and $x$ even, 
or with at least one of $a,b$ even.
\end{theorem}
{\em Remark}. For $n\geq 5, n\neq 6$, the above theorem lefts out the solutions $(a,b,x,y)$ with 
$xab$ odd. These are exactly the {\em exceptional solutions} of the equation \ref{main eqn} in the 
terminology of \cite{Pink}; see also the remark at the end of this paper. 

The proof of Theorem \ref{main result} is given in sections \ref{section n=3}, \ref{section n=4} and
\ref{section ngeq5}, where the cases $n=3$, $n=4$ and $n\geq 5$ are respectively considered. 
Our numerous, crucial computations in section \ref{section n=3} have been done mainly with the aid of
\textsc{Magma} \cite{Bosma},\cite{magma-handbook}; to a less extent we have also been aided by the routines
of \textsc{Pari} ($\mathtt{http://pari.math.u-bordeaux.fr/}$).

Note that since $n\geq 3$, it follows that $n$ is either a multiple of $4$
or a multiple of an odd prime $p$, therefore it suffices to study the equation (\ref{main eqn})
when $n=3,4$ or an odd prime $\geq 5$.
Furthermore, note that if $b=0$, then our equation reduces to the equation $x^{2}+5^{a}=y^{n}$,
which is solved in \cite{Liqun2}.
Also, when $a=0$, the equation (\ref{main eqn}) reduces to $x^{2}+11^{b}=y^{n}$ which is solved
in \cite{CDLPS}. 
\section{ Equation (\ref{main eqn}) with $n=3$} \label{section n=3}
This section is devoted to the proof of the following result.
\begin{proposition} \label{proposition n=3}
The complete solution of the equation
\begin{equation}  \label{n=3,initial}
x^2 +5^a 11^b=y^3\,,\quad a\geq 0,\,b\geq 0,\,x> 0,\,y>0\,, \mathrm{gcd}(x,y)=1 
\end{equation}
is given in (\ref{complete solution}).
\end{proposition}
Writing in (\ref{n=3,initial}) $a=6A+i$, $b=6B+j$ with $0\leq i,j\leq 5$ we see that 
\begin{equation*}
\left(\frac{y}{5^{2A}11^{2B}},\frac{x}{5^{3A}11^{3B}}\right)
\end{equation*}
is an $S$-integral point $(X,Y)$ on the elliptic curve 
\begin{equation*}
{\mathcal{E}}_{ij}: Y^2 = X^3 - 5^i 11^j \:,
\end{equation*}
where $S=\{5,11\}$, with the numerator of $X$ being prime to 55, in view of the restriction $\gcd(x,y)=1$. 
A practical method for the explicit computation of all 
$S$-integral points on a Weierstrass elliptic curve has been developed by
Peth\H{o}, Zimmer, Gebel and Herrmann in \cite{PZGH} and has been
implemented in \textsc{Magma}. The relevant routine $\mathtt{SIntegralPoints}$
worked without problems for all $(i,j)$ except for $(i,j)=(2,5),(4,4),(5,4)$. 
Thus, in the non-exceptional cases $(i,j)$, i.e. 
when $(a,b)\not\equiv (2,5),(4,4),(5,4)\pmod{6}$, all solutions to equation 
(\ref{n=3,initial}) turned out to be those appearing in (\ref{complete solution}). 
For the exceptional pairs $(i,j)=(2,5),(4,4),(5,4)$ \textsc{Magma} returns
no $S$-integral points {\em under the assumption that the rank of the
corresponding curve ${\mathcal{E}}_{ij}$ is zero}, an assumption that the
routine itself cannot certify. Again using \textsc{Magma}, we performed a
2-descent, followed by a 4-descent which proved that the rank is actually
zero in the first two cases $(i,j)=(2,5),(4,4)$, allowing us to arrive
safely to the following conclusion:
\begin{quote}
When $(a,b)\not\equiv (5,4)\pmod{6}$, all solutions to equation (\ref{n=3,initial}) 
are those displayed in (\ref{complete solution}).
\end{quote}
In the third exceptional case $(i,j)=(5,4)$, the 4-descent reveals the
non-torsion point 
\begin{equation*}
(X,Y) = (\frac{997597438498050698749}{101288668233063249}, 
\frac{31508127105495852851671290908932}{32236010714473507582283943})
\end{equation*}
on the curve ${\mathcal{E}}_{54}$, which proves invalid the assumption under
which \textsc{Magma} ``claims'' non-existence of $S$-integral points on 
${\mathcal{E}}_{54}$. Thus, non-existence of integral solutions to 
(\ref{n=3,initial}) when $(a,b)\equiv (5,4)\pmod{6}$ cannot be considered as a fact
that has been proved by \textsc{Magma} routines. Therefore we treat this
equation separately, indicating thus an alternative method for resolving
equations $x^2 + C = y^3$ when $C$ has a prescribed (``small'') set of
distinct prime divisors. Moreover, the resolution of the Thue-Mahler equation (\ref{eq TM2}) that we present 
in section \ref{proposition n=3} is interesting {\em per se}, as it deals with a totality
of non-trivial computational problems that never before (to the best of our knowledge) have been encountered
in the resolution of a Thue-Mahler equation; we acknowledge here the great usefulness of the relevant 
routines of \textsc{Magma}.

In conclusion, according to our discussion so far, for the proof of Proposition \ref{proposition n=3}
it remains to show that the equation 
\begin{equation}  \label{n=3, special}
x^2 +5^a 11^b=y^3\,,\quad (a,b)\equiv (5,4)\!\!\!\pmod{6},\,x\neq 0,\,
\mathrm{gcd}(x,5\cdot 11)=1 \:,
\end{equation}
has no solutions. We write (\ref{n=3, special}) as  
\begin{equation}  \label{n=3, special,modified}
y^3 -5^2\cdot 11\cdot(5^c 11^d)^3 =x^2\,,\quad%
\mbox{$cd$ odd, $(y,5\cdot 11)=1$.}
\end{equation}
and in what follows we will reduce its solution to a number of Thue or Thue-Mahler equations. 
A practical solution of Thue equations has been developed by Tzanakis and de Weger \cite{TdW1} 
which later was improved by Bilu and Hanrot \cite{BH} and implemented in \textsc{Pari} and \textsc{Magma}. 
We will make use of the relevant routines several times without special mentioning. 
Concerning the Thue-Mahler equations, no automatic resolution is available so far and we will follow the method 
of Tzanakis and de Weger \cite{TdW2}.

Factorization of (\ref{n=3, special,modified}) in the field $\Q(\theta)$, 
where $\theta^3=5^2\cdot 11$, gives 
\begin{equation*}
(y-5^c11^d\theta)(y^2+5^c11^dy\theta+5^{2c}11^{2d}\theta^2)=x^2 \:.
\end{equation*}
In the field $\Q(\theta)$ the ideal class-number is 3, an integral
basis is given by $1,\theta,\theta^2/5$ and the fundamental unit is $\epsilon=1+338\theta-52\theta^2$ 
with norm $+1$. It is easily checked that the two factors in the left-hand side of the last equation above 
are relatively prime, hence we have an ideal equation $(y-5^c11^d\theta)=\mathfrak{a}^2$, 
where $\mathfrak{a}$ is an integral ideal. Since the
class-number is relatively prime to the exponent of $\mathfrak{a}$, this
ideal must be principal, generated by an integral element $u+v\theta+w\theta^2/5$. 
Then, passing to element equation, we get 
\begin{equation*}
y-5^c11^d\theta =\pm\epsilon^i (u+v\theta+w\theta^2/5)^2 \:.
\end{equation*}
Taking norms we see that, necessarily, the plus sign must hold above. Also,
comparing coefficients of $\theta$ in both sides we see very easily that $w$
must be divisible by 5, hence, on replacing $w$ by $5w$, we rewrite the last
equation as follows: 
\begin{equation}  \label{element eqn}
y-5^c11^d\theta =\epsilon^i (u+v\theta+w\theta^2)^2 \:.
\end{equation}
We consider two cases, depending on the value of $i$.

Let $i=0$. Equating coefficients of like powers of $\theta$ in both sides of
(\ref{element eqn}) we obtain the following relations:
\begin{eqnarray}
v^2+2uw & = & 0  \label{i=0, coeff2} \\
2uv+275w^2 & = & -5^c\cdot 11^d  \label{i=0, coeff1} \\
u^2+550vw & = & y  \label{i=0, coeff0}
\end{eqnarray}
The above equations, along with the fact that $\gcd(y,5\cdot 11)=1$, easily
imply that $\gcd(u,w)=1$ and $w$ is odd, hence (\ref{i=0, coeff2}) implies
that 
\begin{equation*}
u=2sv_1^2\,,w=-sv_2^2\,,v=2v_1v_2\,,\quad s\in\{-1,1\}\,,\gcd(2v_1,v_2)=1
\end{equation*}
and substitution into (\ref{i=0, coeff1}) gives 
\begin{equation}  \label{factor coeff1}
-5^c\cdot 11^d = 275v_2^4+8sv_1^3v_2= v_2((2sv_1)^3+275v_2^3) \:.
\end{equation}
Since $c$ is odd, it is $\geq 1$, therefore, from (\ref{factor coeff1}) one
of $v_1,v_2$ is divisible by 5. If 5 divides $v_1$, then 25 divides the
right-hand side, hence $c$ must be $\geq 2$ and since it is odd, it must be
at least 3. But then $5^3$ divides $275v_2^4$, hence 5 divides $v_2$ which
contradicts the fact that $\gcd(v_1,v_2)=1$. Therefore, 5 divides $v_2$ and
does not divide $v_1$. \newline
If 11 also divides $v_2$, then, neither 5 nor 11 divides $(2sv_1)^3+275v_2^3$,
 hence, by (\ref{factor coeff1}) we have $v_2=\pm 5^c\cdot 11^d$ and $(2sv_1)^3+275v_2^3=\pm 1$. 
But since the only solutions of the Thue equation 
$X^3+275Y^3=\pm 1$ is $(X,Y)=(\pm 1,0)$, the previous equation is
impossible. 
Therefore 11 divides $v_1$ and does not divide $v_2$. Now, if $d>1$, then 
(\ref{factor coeff1}) implies that $275v_2^4$ is divisible by $11^2$, hence
11 divides $v_2$, a contradiction. Therefore, $d=1$ and by our discussion so
far we conclude, in view of (\ref{factor coeff1}), that $v_2=\pm 5^c$ and 
$(2sv_1)^3+275v_2^3 =\mp 11$, which is impossible since the Thue equation 
$X^3+275Y^3=11$ is impossible.

We conclude therefore that equation (\ref{element eqn}) is impossible when $i=0$.

Next, let $i=1$. Equating coefficients of like powers of $\theta$ in (\ref{element eqn}) gives 
\begin{equation}  \label{i=1,coeff2}
-52u^2+676vu+2wu+v^2+92950w^2-28600wv = 0
\end{equation}
\begin{equation}  \label{i=1,coeff1}
338u^2+2vu-28600wu-14300v^2+275w^2+185900wv = -5^c\cdot 11^d
\end{equation}
\begin{equation}  \label{i=1,coeff0}
u^2-28600vu+185900wu+92950v^2-3932500w^2+550wv = y
\end{equation}
From the above equations it is easy to see that $v$ is even $w$ is odd and
(since also $\gcd(y,5\cdot 11)=1$) $\gcd(u,w)=1$. On the other hand,
equation (\ref{i=1,coeff2}) can be written as 
\begin{equation*}
2(52u-2199 w)(1099u-46475w)=(v+338u-14300w)^2 \:.
\end{equation*}
Since 
\begin{equation*}
\left|\begin{array}{rr}
                 52 & -2199 \\ 
                 1099 & -46475
         \end{array}\right| = 1
\end{equation*}
and $\gcd(u,w)=1$, it follows that the two parenthesis in the left-hand side
of the last equation are relatively prime, the first one being odd, because 
$w$ is odd. It follows that 
\begin{eqnarray*}
52u-2199w & = & sX^2 \\
1099u -46475 w & = & 2sY^2 \\
v+338u-14300w & = & \pm 2XY
\end{eqnarray*}
where $X,Y$ are integers and $s\in\{-1,1\}$. Solving the system in $u,v,w$
we obtain expressions of $u,v,w$ in terms of $X,Y$; then, substitution into 
(\ref{i=1,coeff1}) gives 
\begin{equation*}
150975X^4\pm 185900X^3Y+85800X^2Y^2 \pm 17592XY^3+1352Y^4 = 5^c\cdot 11^d \:.
\end{equation*}
Replacement of $-X$ by $X$ shows that we may consider only the plus sign in
the above equation. We have thus obtained a Thue-Mahler equation which we
will solve in the next section.
%
\subsection{\normalsize The solution of the Thue-Mahler equation} \label{Thue-Mahler}

In this section we prove that the Thue-Mahler equation 
\begin{equation}  \label{eq TM1}
150975X^4 + 185900X^3Y+85800X^2Y^2 + 17592XY^3+1352Y^4 = 5^c\cdot 11^d
\end{equation}
has no solutions.
We will follow closely the method of \cite{TdW2} which, to the best of our knowledge is the only systematic
exposition found so far in the literature. For the convenience of the reader, we will use the same notation 
with \cite{TdW2} as far as possible. {\em The notation in this section is independent of the notation used in the
others sections of the present paper.}

Putting $x=2\cdot 13^2Y, y=X$ (obviously, $(x,y)=1$), we transform equation (\ref{eq TM1}) into 
\begin{equation} \label{eq TM2}
x^4+4398x^3y+7250100x^2y^2+5309489900xy^3+1457454977550y^4 = 2\cdot 13^6 5^c 11^d. 
\end{equation}
We work in the field $K=\Q(\th)$, where $\th$ is a root of the polynomial 
\[
 g(t)= t^4+4398t^3+7250100t^2+5309489900t+1457454977550 \in\Q[t].
\]
The ideal-class number is 1 and an integral basis is $1,\th, (4\th+\th^2)/169$, $(92950\th+173\th^2+\th^3)/142805$.
For shortness, we will use the notation $\ga=[a,b,c,d]$, where $a,b,c,d\in\Z$, to mean that the algebraic integer
$\ga\in K$ has the coefficients $a,b,c,d$ with respect to the above integral basis. 
\\
A pair of fundamental units is 
\begin{eqnarray*}
 \ep_1 & = & [677070473, 1764897, 260182,69044] \\
\ep_2 & = & [7564704083, 22782192, 3852447, 1164105]\,.
\end{eqnarray*}
The factorization of the rational primes 2, 5, 11 and 13 is as follows:
\begin{align*}
2=\ep_2^{-1}\pi_2^4\, & \quad \pi_2=[21436,39,3,0] \\
 5=\ep_1^{-1}\ep_2\pi_{51}\pi_{52}^3\,, & \quad \pi_{51}=[9690469,26053,3965,1087] \\
                 & \quad \pi_{52}=[653350925,1762426,269424,74288] \\
11=\pi_{111}\pi_{112}^3\,, & \quad \pi_{111}=[1060859,2835,429,117] \\
                   & \quad \pi_{112}=[204919,535,79,21] \\
13 =-\pi_{131}\pi_{132}\,, & \quad \pi_{131}=[127759589,344590,52671,14521] \\
                 & \quad \pi_{132}=[16961503,45062,6773,1833]\,,
\end{align*}
where all the prime elements $\pi_{ij}$ above, except for $\pi_{132}$ are of degree 1, and $\pi_{132}$ is
of degree 3.
\\
In the notation of relation (3) of \cite{TdW2}, 
\[
 f_0\leftarrow 1,\,c\leftarrow 2\cdot 13^6,\,p_1\leftarrow 5,\, z_1\leftarrow c,\,
p_2\leftarrow 11, z_2\leftarrow d\,.
\]
Fix the prime $p\in\{5,11\}$. For the elements of the ring $\Z_p$ of $p$-adic integers we will use the 
notation $0.d_0d_1d_2\ldots$, where $d_0,d_1,d_2,\ldots$ are integers between $0$ and $p-1$, to mean 
the $p$-adic integer $d_0+d_1p+d_2p^2+\cdots$. In the sequel, all our computations with $p$-adic numbers
have been done with the relevant \textsc{Magma} routines.
\\
Over $\Q_p[t]$ we have the factorization of $g(t)=g_1(t)g_2(t)$ as in the following table.
\begin{center}
 \begin{supertabular}{|c|c|c|}
\multicolumn{3}{c}{Factorization $g(t)=g_1(t)g_2(t)$ into irreducibles over $\Q_p$} \\[2mm]
 \hline \hline
&                     &            \\[-4mm]
$p$ & $g_1(t)$ & $g_2(t)$ \\  \hline\hline
5 & $t-(0.20404\ldots)$ & $t^3+(0.00011\ldots)t^2+(0.00422\ldots)t+(0.00444\ldots)$ \\ \hline
11 & $t-(0.25033\ldots)$ & $t^3 +(0.09363\ldots)t^2 +(0.09900\ldots)t+(0.052(10)6\ldots)$\\ \hline
\end{supertabular}
\end{center}
We denote by $L_p$ the splitting field of $g(t)$ over $\Q_p$. This is obtained in two steps, 
as shown in the following table.
\begin{center}
\begin{supertabular}{|c|c|c|} 
\multicolumn{3}{c}{$K_p=\Q_p(u)$, $g_{p1}(u)=0$ and $L_p=K_p(v)$, $g_{p2}(v)=0$} \\[2mm]
 \hline \hline
&                     &            \\[-4mm]
$p$ & $g_{p1}(t)$ & $g_{p2}(t)$ \\  \hline\hline
5 & $t^2+4t+2$ & $t^3+(0.03001\ldots)t^2+(0.00111\ldots)t+(0.04442\ldots)$ \\ \hline
11 & $t^2+7t+2$ & $g_2(t)$\\ \hline
\end{supertabular}
\end{center}
The roots of $g(t)$ are shown in the following table.
\begin{center}
 \begin{supertabular}{|c|c|c|}
\multicolumn{3}{c}{The roots of $g(t)$ over $\Q_p$} \\[2mm]
 \hline \hline
&                     &            \\[-4mm]
$p$ & $\cnj{\th}{1}$ & $\cnj{\th}{2},\,\cnj{\th}{3},\,\cnj{\th}{4}$  \\ \hline\hline
    &           &    \\[-2mm]
5 & $0.20404\ldots$ & $(0.10001\ldots)v^2+(0.04220\ldots)v+(0.00011\ldots)$, \\[1mm]
  &                 & $(0.30230\ldots)uv^2+(0.33140\ldots)v^2+(0.00424\ldots)uv$ \\ 
  &                 &  $\hspace{4mm}+(0.03143\ldots)v+(0.00344\ldots)u+(0.00132\ldots)$, \\[1mm]
  &                & $(0.24214\ldots)uv^2+(0.11303\ldots)v^2+(0.00120\ldots)uv$ \\
  &                &    $\hspace{4mm}+(0.03032\ldots)v+(0.00203\ldots)u+(0.00444\ldots)$ \\ \hline
11 & $0.25033\ldots$ & $v$, \\[1mm]
  &                 & $(0.(10)4(10)71\ldots)uv^2+(0.26306\ldots)v^2+(0.(10)3(10)63\ldots)uv$ \\ 
  &                 &  $\hspace{4mm}+(0.72327\ldots)v+(0.026(10)3\ldots)u+(0.08801\ldots)$, \\[1mm]
  &                & $(0.16039\ldots)uv^2+(0.947(10)4\ldots)v^2+(0.17047\ldots)uv$ \\
  &                &    $\hspace{4mm}+(0.38783\ldots)v+(0.09407\ldots)u+(0.05936\ldots)$ \\ \hline
 \end{supertabular}
\end{center}

For $\ga \in \Q_p$ we define, as usually, $\ordp(\ga)=m$ iff $\ga=p^m\mu$, where $\mu$ is a $p$-adic unit.
We extend the function $\ordp$ to $L_p$ by the formula
\[
 \ordp(\ga)=\frac{1}{6}\ordp(\Nm{L_p}{\Q_p}(\ga))
\]
(see section 4 of \cite{TdW2}). By Statement (i) of the  {\em First Corollary of Lemma 1}, p.231 of \cite{TdW2} 
we conclude that at most one among $\pi_{51}$ and $\pi_{52}$ divides $x-y\th$. Moreover, if $\pi_{52}$ divides
$x-y\th$, then Statement (ii) of the same {\em Corollary} asserts that at most $\pi_{52}^2$ divides $x-y\th$.
Similarly, at most one among $\pi_{111}$ and $\pi_{112}$ divides $x-y\th$ and if this is the case with the 
second one, then at most its first power divides $x-y\th$. 
These observations and standard arguments of Algebraic Number Theory lead to the following 
ideal equation (cf. relation (9) of \cite{TdW2})
\[
 \ideal{x-y\th} =\ida\idb\idp_1^{n_1}\idp_2^{n_2}\,,
\]
where
\[
\ida\in\{\ideal{\pi_2\pi_{131}^6},\ideal{\pi_2\pi_{131}^3\pi_{132}},\ideal{\pi_2\pi_{132}^2}\},
\]
\[
 \idb=\ideal{\pi_{52}^{j_1}\pi_{112}^{j_2}},\quad (0\leq j_1\leq 2,\,0\leq j_2\leq 1)
\]
\[
\idp_1=\ideal{\pi_{51}},\quad  \idp_2=\ideal{\pi_{52}}\,.
\]
\begin{eqnarray} 
 c=n_1+j_1 & \mbox{with ($n_1>0$ and $j_1=0$) or ($n_1=0$ and $j_1\leq 2$),} \nonumber \\[-1mm]
                &                                        \label{c,d} \\[-1mm]
 d=n_2+j_2 & \mbox{with ($n_2>0$ and $j_2=0$) or ($n_2=0$ and $j_2\leq 1$).} \nonumber
\end{eqnarray}
Following the strategy (and notation) of section 7 of \cite{TdW2} we obtain the following relation:
\begin{eqnarray} 
\la & = & \de_2\left(\frac{\cnj{\ep_1}{i_0}}{\cnj{\ep_1}{j}}\right)^{a_1} 
           \left(\frac{\cnj{\ep_2}{i_0}}{\cnj{\ep_2}{j}}\right)^{a_2} 
           \left(\frac{\cnj{\pi_{51}}{i_0}}{\cnj{\pi_{51}}{j}}\right)^{n_1}
           \left(\frac{\cnj{\pi_{111}}{i_0}}{\cnj{\pi_{112}}{j}}\right)^{n_2} \nonumber \\[-1mm]
        &                                                    \label{eq lambda} \\[-1mm]
 \la & = & \de_1\left(\frac{\cnj{\ep_1}{k}}{\cnj{\ep_1}{j}}\right)^{a_1} 
           \left(\frac{\cnj{\ep_2}{k}}{\cnj{\ep_2}{j}}\right)^{a_2} 
           \left(\frac{\cnj{\pi_{51}}{k}}{\cnj{\pi_{51}}{j}}\right)^{n_1}
           \left(\frac{\cnj{\pi_{111}}{k}}{\cnj{\pi_{112}}{j}}\right)^{n_2} - 1\,, \nonumber
\end{eqnarray}
where
\[
 \de_1=\frac{\cnj{\th}{i_0}-\cnj{\th}{j}}{\cnj{\th}{i_0}-\cnj{\th}{k}}\cdot\frac{\cnj{\al}{k}}{\cnj{\al}{j}}\,,
\quad
\de_2=\frac{\cnj{\th}{j}-\cnj{\th}{k}}{\cnj{\th}{k}-\cnj{\th}{i_0}}\cdot\frac{\cnj{\al}{i_0}}{\cnj{\al}{j}}\,,
\]
\[
 \al\in\{\pi_{131}^6,\pi_{131}^3\pi_{132},\pi_{132}^2\}\cdot\pi_2\pi_{52}^{j_1}\pi_{112}^{j_2}
\]
with $j_1,j_2$ as in (\ref{c,d}).
\\
We view (\ref{eq lambda}) either as a relation in $L_p$, where $p\in\{5,11\}$ as the case may be, and in this
case $i_0=1, j=2,k=3$ (cf. table of roots of $g(t)$ over $\Q_p$), or as a relation in $\C$, in which case 
we number the real/complex roots of $g(t)$ as $\th_1,\th_2$, the real ones, and $\th_3,\th_4$ the pair of 
complex-conjugate roots, and we take $i_0\in\{1,2\}$ (one must consider both cases), $j=3,k=4$.
\\
We set now
\[
 A=\max\{|a_1|,|a_2|\},\quad N=\max\{n_1,n_2\},\quad H=\max\{A,N\}.
\]
We apply Yu's theorem (Theorem 1 in \cite{Yu}) in its somewhat simplified version 
presented in Appendix A2 of \cite{TdW2}. Note that the least field generated over $\Q$ by the five algebraic 
numbers appearing in $\la$ (cf. (\ref{eq lambda})) is of degree 24. Thus, in the notation of the above mentioned
Appendix A2, $n_1\leftarrow 24, n_2\leftarrow 48, (q,u)\leftarrow (2,2), f_2\leftarrow 1$, and after the
computation of the various parameters in Yu's theorem, we obtain the values $c_{13},c_{14}$ (p. 238 of \cite{TdW2})
for which
\[
 N \leq c_{13}(\log H +c_{14})\,,
\]
namely, $c_{13}=2.0564\times 10^{39}, c_{14}=7.7425$.
\\
Next, we work with real/complex linear forms in logarithms of agebraic numbers. In the terminology of 
\cite{TdW2}, p. 243, we encounter a ``complex case'', therefore, following that paper, we consider the linear
form 
\[
\Lambda_0=i^{-1}\Log(1+\la)=i^{-1}
           \Log\left(\de_1\left(\frac{\cnj{\ep_1}{k}}{\cnj{\ep_1}{j}}\right)^{a_1} 
           \left(\frac{\cnj{\ep_2}{k}}{\cnj{\ep_2}{j}}\right)^{a_2} 
           \left(\frac{\cnj{\pi_{51}}{k}}{\cnj{\pi_{51}}{j}}\right)^{n_1}
           \left(\frac{\cnj{\pi_{111}}{k}}{\cnj{\pi_{112}}{j}}\right)^{n_2}\right) 
\]
where $\Log$ denotes the principal branch of the complex logarithmic function. Since, for every 
$z\in\C$, $i^{-1}\Log(\overline{z}/z)=\Arg(\overline{z}/z)$, where $\Arg$ denotes the principal Argument,
we have after expansion (remember that $i_0=1,2$ and $j=3,k=4$),
\begin{align}
 \Lambda_0 = 
\Arg\frac{\cnj{\th}{i_0}-\cnj{\th}{3}}{\cnj{\th}{i_0}-\cnj{\th}{4}}\cdot\frac{\cnj{\al}{4}}{\cnj{\al}{3}}
 & +a_1\Log(\de_1\left(\frac{\cnj{\ep_1}{4}}{\cnj{\ep_1}{3}}\right) 
+a_2\Arg \left(\frac{\cnj{\ep_2}{4}}{\cnj{\ep_2}{3}}\right)   \nonumber \\
 & \; +n_1\Arg\left(\frac{\cnj{\pi_{51}}{4}}{\cnj{\pi_{51}}{3}}\right)
  +n_2\Arg\left(\frac{\cnj{\pi_{111}}{k}}{\cnj{\pi_{112}}{j}}\right) +a_0(2\pi)\,.  \label{Lambda0}
\end{align}
According to relation (27), p. 245 of \cite{TdW2}, we have
\[
 0< |\Lambda_0| < 1.02c_{21}e^{-c_{15}A} \,,
\]
where $c_{21}$ and $c_{16}$ are explicit, and we need further a lower bound of the shape 
$|\Lambda_0|>\exp(-c_7(\log H+2.5))$ (see p. 246 of \cite{TdW2}). Baker-W\"{u}stholz's theorem \cite{BW}
furnishes us $c_7=8.43\times 10^{56}$. The constant $c_{16}$ must be less than $3.809\ldots$ (for the
choice of $c_{16}$ see \cite{TdW2}, bottom of p. 239). The constants $c_7,c_{13},c_{14}$ and $c_{16}$
are the crucial ones in the computation of an upper bound for $H$. A number of other parameters must be 
obtained by elementary but quite cumbersome calculations. A very detailed exposition of how this list of parameters
are calculated for the general Thue-Mahler equation is exposed in the first eleven sections of \cite{TdW2},
culminating to an explicit upper bound for $H$. Fortunately, the computation of these constants, including
that of $c_7,c_{13},c_{14}$ and $c_{16}$, can be rather easily implemented in (for example) {\sc Maple} or
{\sc Magma}. For our equation the upper bound that we calculate is
\begin{equation}  \label{K0}
 H < K_0=5.792\times 10^{58}\,.
\end{equation}
Especially for $N$, a smaller upper bound is obtained by the {\em Corollary to Theorem 10}, p. 248 of \cite{TdW2},
namely
\begin{equation}  \label{N0}
 N < N_0 = 2.942\times 10^{41}\,.
\end{equation}
Thus, the upper bound for $H$ is, actually, the upper bound for $A$.

{\bf Reduction of the upper bound}. In order to considerably reduce the upper bound (\ref{N0}) by the so
called {\em $p$-adic reduction process}, we need the {\em $p$-adic logarithmic function} $\log_p z$, which
is defined for every $p$-adic unit $z\in L_p$ and takes values in $L_p$; see the detailed exposition in 
section 12 of \cite{TdW2}.
\\
For $p\in\{5,11\}$ we put (viewing $\la$ in (\ref{eq lambda}) as an element of $L_p$)
\[
 \Lambda =\log_p(1+\la)= \log_p\de_1 +n_1\log_p\frac{\cnj{\pi_{51}}{k}}{\cnj{\pi_{51}}{j}}
                                     +n_2\log_p\frac{\cnj{\pi_{111}}{k}}{\cnj{\pi_{111}}{j}} 
                                     +a_1\log_p\frac{\cnj{\ep_1}{k}}{\cnj{\ep_1}{j}}
                                     +a_2\log_p\frac{\cnj{\ep_2}{k}}{\cnj{\ep_2}{j}}\,,
\]
where the indices can be chosen arbitrarily from the set $\{2,3,4\}$ ($k\neq j$). Expressing $\Lambda$ in
terms of the basis $1,u,v,uv,v^2,uv^2$ of $L_p/\Q_p$, we can write
\[
 \Lambda = \Lambda_0 + \Lambda_1u + \Lambda_2v + \Lambda_3uv + \Lambda_4v^2 +\Lambda_5uv^2\,,
\]
where each $\Lambda_i$ is a linear form
\[
 \Lambda_i=\al_{i0}+\al_{i1}n_1 + \al_{i2}n_2+\al_{i3}a_1 +\al_{i4}a_2\,,\quad \al_{ij}\in\Z_p\,,
\quad (i=0,\ldots,5).
\]
Following the discussion of section 14 of \cite{TdW2}, for each $i$, we divide by the $\al_{ij}$ whose $\ordp$
has a minimal value (actually, this is obtained for som $j>0$), obtaining thus a linear form
\[
 \Lambda'_i =-\be_0-\be_1b_1-\be_2+b_2-\be_3b_3+b_4\,,
\]
where $(b_1,b_2,b_3,b_4)$ is a permutation of $(n_1,n_2,a_1,a_2)$. At this point we note that the $p$-adic numbers
$\be_0,\ldots,\be_3$ are computed with a high $p$-adic precision $m$. We denote by $\cnj{\be}{m}$ the 
{\em rational integer} which approximates $\be$ with $m$ $p$-adic digits; in other words,
$\ordp(\be-\cnj{\be}{m})\geq m$.
\\
Following the $p$-adic reduction process described in section 15 of \cite{TdW2}, we consider the lattice whose 
basis is formed by the columns of the matrix
\[
 \left(\begin{array}{cccc}
        W & 0 & 0 & 0 \\
        0 & 1 & 0 & 0 \\
        0 & 0 & 1 & 0 \\
      \cnj{\be_1}{m} & \cnj{\be_2}{m} & \cnj{\be_3}{m} & p^m
       \end{array}\right)\,,
\]
where $W$ is an integer somewhat larger than $K_0/N_0$; in this case we choose $W=2\cdot 10^{17}$.
Then we obtain an {\em LLL-reduced basis} of the lattice. As explained in section 15 of \cite{TdW2}, if
$m$ is sufficiently large, then it is highly probable that a certain condition stated in Proposition 15
of \cite{TdW2} (in which condition the reduced basis is, of course, involved) is fulfilled; and if the
condition is fulfilled, then, according to that Proposition 15, $n_1,n_2\leq m+1$. It turns out that, if
$p=5$, then $m=306$ is sufficient for the condition of Proposition 15 to be fulfilled; and if $p=11$,
it suffices to have a precision of $m=207$ 11-adic digits. Thus, in the first $p$-adic step we made a huge
``jump'', falling from (\ref{N0}) to $N\leq N_1=307$.
\\
Now it is the turn of the {\em real reduction step}. We rewrite the linear form $\Lambda_0$ in (\ref{Lambda0})
as
\[
\Lambda_0=\rho_0+n_1\la_1+n_2\la_2+a_1\mu_1+a_2\mu_2+a_0\mu_3 \quad (\mu_3=2\pi)
\]
and for $C=10^m$, with $m$ a sufficiently large integer (having nothing to do with the $m$ in the 
$p$-adic reduction process), we put 
\[
 \phi_0=[C\rho_0],\,\phi_i=[C\la_i]\; (i=1,2)\,, \psi_i=[C\mu_i]\; (i=1,2,3)\,,
\]
where $[x]=\lfloor x\rfloor$ if $x\geq 0$ and $[x]=\lceil x\rceil$ if $xN 0$. In practice, this means that we
must compute our real numbers $\rho,\la,\mu$ with a precision of somewhat more than $m$ decimal digits.
\\
Following the discussion of section 16 of \cite{TdW2}, we consider the lattice whose basis is formed by the
columns of the matrix
\[
 \left(\begin{array}{ccccc}
 W & 0 & 0 & 0 & 0 \\
 0 & W & 0 & 0 & 0  \\
 0 & 0 & 1 & 0 & 0 \\
 \phi_1 & \phi_2 & \psi_1 & \psi_2 & \psi_3
\end{array} \right)\,.
\]
Again, we compute an LLL reduced basis for the lattice and, according to Proposition 16 of \cite{TdW2}, if a
certain condition, in which the reduced basis is involved, is satisfied, then a considerably smaller upper bound
for $H$ is obtained. It is highly probable that this condition is satisfied if $C=10^m$ is sufficiently large.
As it turns out in our case, $m=200$ is sufficient and the reduced upper bound implied by the above
mentioned Proposition 16 is $H\leq K_1=546$, an enormous ``jump'' from (\ref{K0})!

This strategy of a $p$-adic reduction process followed by a real reduction process is repeated, with
$K_1$ in place of $K_0$ and $N_1$ in place of $N_0$, giving even smaller upper bounds, namely,
$N\leq N_2=32$ and $H\leq K_2=74$. We repeat the process once more. The 5-adic reduction process gives
$n_1\leq 25$ and the 11-adic reduction process gives $n_2\leq 18$. The real reduction process gives
$H\leq 59$. Thus,
\begin{equation} \label{final bounds}
 0\leq n_1\leq 25\,,\quad 0\leq n_2\leq 18\,,\quad A=\max\{|a_1|,|a_2|\}\leq 59\,.
\end{equation}

{\bf The sieve after the reduction}. The bounds (\ref{final bounds}) cannot be further improved, therefore we have
to search whether there exist quadruples $(a_1,a_2,n_1,n_2)$ in the range (\ref{final bounds}) and pairs
$(i_1,i_2)$ and $(j_1,j_2)$ in the range
\begin{equation} \label{range of i and j}
(i_1,i_2)\in\{(6,0),(3,1),(0,2)\}\,,\quad 0\leq j_1\leq 2\,,\; 0\leq j_2\leq 1\,, 
\end{equation}
such that 
\begin{equation} \label{h}
h(i_1,i_2,j_1,j_2,a_1,a_2,n_1,n_2):=\pi_2\pi_{131}^{i_1}\pi_{132}^{i_2}\pi_{52}^{j_1}\pi_{112}^{j_2}
                                    \ep_1^{a_1}\ep_2^{a_2}\pi_{51}^{n_1}\pi_{111}^{n_2}
\end{equation}
is of the form $x-y\th$, i.e., such that, after expanding the right-hand side in (\ref{h}) and expressing
it in terms of the basis $1,\th,\th^2,\th^3$ of $K/\Q$, the coefficients of $\th^2$ and $\th^3$ are zero.
Doing this check by ``brute force'' is very time consuming. Instead, we choose to do the following
sieving process (see also section 18 of \cite{TdW2}).

Let $q$ be a rational prime which splits into four distinct (first degree) prime divisors $\rho_1,\ldots,\rho_4$
of $K$. Then, for every algebraic integer $\ga\in K$, there exist rational integers $A_i,\,i=1,\ldots 4$, 
such that $\ga\equiv A_i\pmod{\rho_i}$. As a consequence, every (rational) relation with algebraic integers
of $K$ implies congruences $\bmod{\,\rho_i}$, one for every $i=1,\ldots,4$. But since in these congruences
the elements of $K$ are replaced by {\em rational integers}, these are valid also as congruences in $\Z$
modulo $q$.
\\
Take, for example, $q=31$. Then, we have the ideal factorization $\ideal{q}=\prod_{i=1}^4\ideal{\rho_i}$,
where
\[
 \th\equiv 1\!\!\!\pmod{\rho_1}\,,\; \th\equiv 17\!\!\!\pmod{\rho_2}\,,\; 
\th\equiv 19\!\!\!\pmod{\rho_3}\,,\;\th\equiv 29\!\!\!\pmod{\rho_4}\,.
\]
From a relation of the form $x-y\th=h(\bi)$, where $\bi=(i_1,i_2,j_1,j_2,a_1,a_2,n_1,n_2)$ is in
the range (\ref{final bounds}) and (\ref{range of i and j}), we obtain the four congruences
\begin{equation}\label{H congruences}
 x-y\equiv H_1(\bi)\,,\; x-17y\equiv H_2(\bi)\,,\; x-19y\equiv H_3(\bi)\,,\; x-29y\equiv H_4(\bi)
\!\!\!\pmod{31}
\end{equation}
where $H_1(\bi)$ is the rational integer resulting on replacing $\th$ by 1 in $h(\bi)$, and similarly
for the remaining $H_j(\bi)$'s. Then,
\begin{equation} \label{relation with H}
 27H_1(\bi)+5H_2(\bi)\equiv H_3(\bi)\,,\quad 7H_1(\bi)+25H_2(\bi)\equiv H_4(\bi)\!\!\!\pmod{31}\,.
\end{equation}
Note now that, for every algebraic integer $\ga\in K$, the order of $\ga$ modulo 31 is a divisor of 30.
The orders of $\ep_1,\ep_2,\pi_{51},\pi_{111}$ modulo 31 are 30,15,15,30, respectively. Therefore, we check
the congruences (\ref{relation with H}) for all $\bi$'s with $(i_1,i_2)$ and $(j_1,j_2)$ as in 
(\ref{range of i and j}) and $0\leq a_1\leq 29$, $0\leq a_2\leq 14$, $0\leq n_1\leq 14$, $0\leq n_2\leq 18$.
\\
For example, when $(i_1,i_2,j_1,j_2)=(6,0,2,1)$, there are 4275 quadruples $(a_1,a_2,n_1,n_2)$
that satisfy the first congruence (\ref{relation with H}). We check which of them also satisfy the second
congruence (\ref{H congruences}) and only 117 quadruples pass the test. Now, these 117 quadruples must be lifted to 
cover the range (\ref{final bounds}), resulting  to 6532 quadruples. Thus, there are 6532 6-tuples
$\bi=(6,0,2,1,a_1,a_2,n_1,n_2)$ with $a_1,a_2,n_1,n_2$ in the range indicated by (\ref{final bounds}), that
satisfy both congruences (\ref{H congruences}).
\\
Next, we work similarly with the prime $q=79$ (which splits into four distinct prime divisors of $K$).
The analogous to the congruences (\ref{H congruences}) are now
\[
x-6y\equiv H_1'(\bi)\,,\; x-14y\equiv H_2'(\bi)\,,\; x-41y\equiv H_3'(\bi)\,,\; x-44y\equiv H_4'(\bi)
\!\!\!\pmod{73}\,,
\]
implying the anlogous to (\ref{relation with H}) congruences
\[
24H_1'(\bi)-23H_2'(\bi)\equiv H_3'(\bi)\,,\quad -22H_1'(\bi)+23H_2'(\bi)\equiv H_4'(\bi) \pmod{73}\,.
\]
We check which of the 6532 6-tuples $\bi$, obtained before, satisfy the last congruences and only three 6-tuples
pass the test, which are tested by a final similar test with the prime 223 in place of 73; no one passes the
test. This shows that no 6-tuple $(6,0,2,1,a_1,a_2,n_1,n_2)$ is accepted.
\\
A similar test is repeated for every $(i_1,i_2)$ and $(j_1,j_2)$ as in (\ref{range of i and j}) and always we
end up with no acceptable $(i_1,i_2,j_1,j_2,a_1,a_2,n_1,n_2)$.

{\bf Final conclusion of section \ref{Thue-Mahler}}: 
The Thue-Mahler equation (\ref{eq TM2}) has no solutions and, consequently, neither  
the equation (\ref{eq TM1}) has solutions. This completes the proof of Proposition \ref{proposition n=3}.
%
%
%
\section{Equation (\ref{main eqn}) with $n=4$} \label{section n=4}
In this section we prove the following result.
\begin{proposition}  \label{proposition n=4}
If $n=4$, then the equation (\ref{main eqn}) has no solution.
\end{proposition}
\proof Since $n=4$, equation (\ref{main eqn}) is written as
\begin{equation} \label{eq with n=4}
5^{a}\cdot 11^{b}=(y^{2}+x)(y^{2}-x)\,,
\end{equation}
from which we obtain
\begin{eqnarray*}
y^{2}+x &=& 5^{a_{1}}11^{b_{1}} \\
y^{2}-x &=& 5^{a_{2}}11^{b_{2}}
\end{eqnarray*}
where $a_{1},a_{2},b_{1},b_{2}\geq 0$.  
From the equations above and the assumption $\gcd(x,y)=1$ it follows that $a_1,a_2$ cannot both be positive, and 
similarly for $b_1,b_2$. Summing the two equations we obtain
\begin{equation} \label{eq de Weger}
Z^{2}-Du^{2}=2\cdot 5^{a_{1}}11^{b_{1}}
\end{equation}
where $D\in\{2,10,22,110\}$, $Z=2y$ and $u=5^{a_{2}}11^{b_{2}}$. 
We have $\gcd(Z,u)=1$. Indeed, otherwise we would have $\gcd(2y^2,u)>1$, hence
$\gcd(5^{a_1}11^{b_1}+5^{a_2}11^{b_2},5^{a_2}11^{b_2})>1$ which contradicts our remark concerning the pairs
$a_1,a_2$ and $b_1,b_2$ a few lines above.
\\
We claim that $a_1=0$. Indeed, suppose that $a_1\geq 1$. If $D=2$ or $22$, then by (\ref{eq de Weger}),
$Z^2\equiv 2u^2\pmod{5}$, implying $Z\equiv 0\equiv 0\pmod{5}$, which contradicts $\gcd(Z,u)=1$. 
If $D=10$ or $110$, then (\ref{eq de Weger}) implies $Z\equiv 0\pmod{5}$, hence also $2y^2\equiv 0\pmod{5}$.
But $2y^2=5^{a_1}11^{b_1}+5^{a_2}11^{b_2}$ and we have assumed that $a_1\geq 1$; therefore, $a_2\geq 1$,
contradicting our remark that $a_1,a_2$ cannot both be positive.
\\
With completely analogous arguments we prove that $b_1=0$, by distinguishing the cases $D=2,10$ and
$D=22,110$ and taking into account that $b_1,b_2$ cannot both be positive.

Thus, $y^2+x=1$, which is impossible since $x$ and $y$ are positive integers.
\proofend

%
\section{Equation (\ref{main eqn}) with $n\geq 5, n\neq 6$} \label{section ngeq5}
In this section we prove the following result.
\begin{proposition}  \label{proposition ngeq5}
The equation 
\begin{equation} \label{eq ngeq5}
x^2+ 5^a 11^b =y^n\,,\quad (x,y)=1,\,n\geq 5 \,,
\end{equation}
is impossible if at least one among $a$ and $b$ is even or if $ab$ is odd and $x$ is even.
\end{proposition}
\proof 
Since in the previous sections we have completed the study of the equation $x^2+ 5^a 11^b =y^n$ with $n=3,4$,
we certainly can assume that $n$ is a prime $\geq 5$.
\\
We write (\ref{eq ngeq5}) as 
\begin{equation} \label{eq with z}
x^2+dz^2 = y^n\,,\quad d\in\{1,5,11,55\}\,,\;z=5^{\al}11^{\be}\,, 
\end{equation}
where the relation of $\al$ and $\be$ with $a$ and $b$, respectively, is clear.
\\
If at least one among $a$ and $b$ is even, then $d\in\{1,5,11\}$ and we see $\bmod{\,8}$ that $x$ is even.
If both $a$ and $b$ are odd, then $d=55$ and both cases, $x$ even or odd can arrise. According to the announcement 
of the Proposition, we consider only the case that $x$ is even.

We work in the field $\Q(\sqrt{-d})$. The algebraic integers in this number field are of the form
$ (u+v\sqrt{-d})/2$, where $u,v\in\Z$ with $u,v$ both even, if $d=1,5$ and 
$u\equiv v\!\!\!\pmod{2}$ if $d=11,55$.
\\
Since $x$ is even, the factors in the left-hand side of the equation $(x+z\sqrt{-d})(x-z\sqrt{-d})=y^n$
are relatively prime and we obtain the ideal equation $\ideal{x+y\sqrt{-d}}=\ida^n$. Then, since the ideal-class
number is 1, 2, or 4, and $n$ is odd, we conclude that the ideal $\ida$ is principal. Moreover, the units are
$\pm 1$ and, in case $d=1$, also $\pm i$ ($i=\sqrt{-1}$). In any case, the units are always $n$-th powers, so
that we can finally write
\[
 x+z\sqrt{-d}=\mu^n\,,\quad \mu=\frac{u+v\sqrt{-d}}{2}\,,
\]
where $u,v\in\Z$, with $u,v$ both even, if $d=1,5$ and $u\equiv v\pmod{2}$ if $d=11,55$. For any 
$\ga\in\Q(\sqrt{-d})$ we denote by $\overline{\ga}$ the conjugate of $\ga$. Note that
\[
 \mu-\overline{\mu}=v\sqrt{-d}\,,\quad \mu+\overline{\mu}=u\,,\quad \mu\overline{\mu}=\frac{u^2+dv^2}{4}\,.
\]
We thus obtain
\begin{equation} \label{eq Lucas}
\frac{2\cdot 5^{\al}11^{\be}}{v}=
 \frac{2z}{v} = \frac{\mu^n-\overline{\mu}^n}{\mu-\overline{\mu}}=\mbox{(by definition) $L_n\in\Z$} \,.
\end{equation}
Thus, $\frac{2z}{v}$ is the $n$-th term of {\em Lucas sequence} $(L_m)_{m\geq 0}$. Note that
\begin{equation} \label{recurence}
 L_0=0,\,L_1=1\,,\; L_m=uL_{m-1}-\frac{u^2+dv^2}{4}L_{m-2}\,,\: m\geq 2\,.
\end{equation}
Following the nowadays standard strategy based on the important paper \cite{PrimDiv}, we distinguish two
cases according as $L_n$ has or has not primitive divisors.
\\
Suppose first that $L_n$ has a primitive divisor, say $q$. 
By definition, this means that the prime $q$ divides $L_n$ and $q$ does not divide 
$(\mu-\overline{\mu})^2L_1\cdots L_{n-1}$, hence 
\begin{equation} \label{q does not divide}
 q\not|(\mu-\overline{\mu})^2L_1\cdots L_4=(dv^2)\cdot u\cdot\frac{3u^2-dv^2}{4}\cdot\frac{u^2-dv^2}{2}\,.
\end{equation}
If $q=2$, then (\ref{q does not divide}) implies that $uv$ is odd, hence $d=11$ or $55$. If $d=11$, then
the third factor in the right-most side of (\ref{q does not divide}) is even, a contradiction. 
If $d=55$, then, from (\ref{recurence}) we see that $L_m\equiv L_{m-1}\pmod{2}$, hence $L_m$ is odd for every 
$m\geq 1$, implying that 2 cannot be a primitive divisor of $L_n$.
\\
If $q=5$, then (\ref{q does not divide}) implies that $d=1,11$ and 5 does not divide
$uv(3u^2-dv^2)(u^2-dv^2)$. It follows easily then that $v^2\equiv-u^2\pmod{5}$, so that, by (\ref{recurence}),
$L_m\equiv uL_{m-1}\pmod{5}$ for every $m\geq 2$. Therefore, $5\not| L_n$, so that 5 cannot be a primitive
divisor of $L_n$.
\\
If $q=11$, then, by (\ref{q does not divide}), $d=1,5$ and we write $u=2u_1,v=2v_1$ with $u_1,v_1\in\Z$, so
that $\mu=u_1+v_1\sqrt{-d}$ and (\ref{q does not divide}) becomes
$q\not|u_1v_1(3u_1^2-dv_1^2)(u_1^2-dv_1^2)$. Moreover, $L_m=2u_1L_{m-1}-(u_1^2+dv_1^2)L_{m-2}$ for $m\geq 2$.
Note that $\mu\overline{\mu}=u_1^2+dv_1^2\not\equiv 0\pmod{11}$; therefore, by Corollary 2.2 of \cite{PrimDiv},
there exists a positive integer $m_{11}$ such that $11|L_{m_{11}}$ and $m_{11}|m$ for every $m$ such that
$11|L_m$. It follows then that $11|\gcd(L_n,L_{m_{11}})=L_{\gcd(n,m_{11})}$\footnote{By the well-known property
of Lucas sequences: $\gcd(L_m,L_k)=L_{\gcd(m,k)}$.}. Because of the minimality property of $m_{11}$, we 
conclude that $\gcd(n,m_{11})$, hence, since $n$ is prime, $m_{11}=n$. On the other hand, the Legendre symbol
$\left(\frac{(\mu-\overline{\mu})^2}{11}\right)=-1$, hence, by Theorem XII of \cite{Carm}
(or by Theorem 2.2.4\,(iv) of \cite{LucaReport}), $11|L_{12}$. Therefore $m_{11}|12$, i.e. $n|12$, a contradiction,
since $n$ is a prime $\geq 5$.
\\
We therefore conclude that $L_n$ has no primitive divisors. Then, by Theorem 1.4 of \cite{PrimDiv}, $n<30$. 
By (\ref{eq Lucas}), the prime divisors of $L_n$ belong to $\{2,5,11\}$ and now, looking at the table 1
of \cite{PrimDiv}, we se that the only possibility is $n=5$ and $(u,-dv^2)=(1,-11)$, i.e. 
$\mu=(1+\sqrt{-11})/2$. Going back to (\ref{eq Lucas}) we obtain no solution.
\proofend

{\em A remark on the case when in (\ref{eq ngeq5})  $a,b$ and $x$ are odd.} We explain here why the method applied
for the proof of Proposition \ref{proposition ngeq5} does not apply when $abx$ is odd. In this case $d=55$ 
and we work in the field $\Q(\th)$, where $\th^2-\th+14=0$. 
The equation (\ref{eq ngeq5}) is factorized as $(x-z+2z\th)(x+z-2z\th)=y^n$, where the factors in the left-hand 
side are {\em not} relatively prime. Then, using rather standard arguments of algebraic Number Theory, we are led to
the equation
\[
\frac{5^{\al}11^{\be}}{v}=z=\frac{1}{2^{2(n+1)}}\cdot
\frac{(1+\th)^{\frac{n+1}{2}}\mu^n-(2-\th)^{\frac{n+1}{2}}\overline{\mu}^n}{\mu-\overline{\mu}}\,,
\]
which is the analogous to equation (\ref{eq Lucas}). Now, however, although the right-hand side is a term
of a second order recurrence sequence, it is not a term of a {\em Lucas sequence} and consequently we cannot
argue based on the results of \cite{PrimDiv} as we previously did.

\subsection*{Acknowledgements}
The first named author was supported by the research fund of Uludag University 
Project No F-2008/31.


\author{Ismail Naci Cang\"{u}l, Musa Demirci \\[-1mm] 
Department of Mathematics, Uluda\u{g} University \\[-1mm] 
16059 Bursa, TURKEY \\[-2mm] 
e-mail: cangul@uludag.edu.tr, mdemirci@uludag.edu.tr
\\[2mm]

\and
G\"{o}khan Soydan \\[-1mm] 
Isiklar Air Force High School, \\[-1mm] 
16039 Bursa, TURKEY \\[-2mm]
e-mail: gsoydan@uludag.edu.tr
\\[2mm]

\and
Nikos Tzanakis \\[-1mm]
Department of Mathematics, University of Crete \\[-1mm] 
71409 Iraklion-Crete, GREECE \\[-2mm] 
e-mail: tzanakis@math.uoc.gr
}

\begin{thebibliography}{99}
\normalsize
\baselineskip=17pt
%
\bibitem{Arif4} S.A.~Arif and F.S.~Abu Muriefah, {\em On Diophantine equation $x^2+2^k=y^n$},  
Int.~J.~Math.~Math.~Sci. {\bf 20}, No 2 (1997), 299-304.
%
\bibitem{Arif1} S.A.~Arif and F.S.~Abu Muriefah, {\em The Diophantine equation $x^2+3^m=y^n$},  
Int.~J.~Math.~Math.~Sci. {\bf 21} (1998), 619-620.
%
\bibitem{Arif2} S.A.~Arif and F.S.~Abu Muriefah, {\em The Diophantine equation $x^2+5^{2k+1}=y^n$},  
Indian J.~Pure Appl.~Math. {\bf 30} (1999), 229-231.
%
\bibitem{Arif3} S.A.~Arif and F.S.~Abu Muriefah, {\em On the Diophantine equation $x^2+q^{2k+1}=y^n$},
J.~Number Th. {\bf 95} (2002), 95-100.
%
\bibitem{FSAbu7} F.S.~Abu Muriefah, {\em On the Diophantine equation $x^2+5^{2k}=y^n$}, 
Demonstratio Math. {\bf 319} No 2 (2006), 285-289.
%
\bibitem{FSAbu1} F.S.~Abu Muriefah and Y.~Bugeaud, {\em The Diophantine equation $x^2+C=y^n$: a brief overview}, 
Rev.~Colombiana Math. {\bf 40} (2006), 31-37.
%
\bibitem{FSAbu3} F.S.~Abu Muriefah, F.~Luca and A.~Togbe, {\em On the Diophantine equation 
$x^{2}+5^{a}13^{b}=y^{n}$}, Glasgow Math.~J. {\bf 50} (2008), 175-181.
%
\bibitem{BW} { A.~Baker and G.~W\"{u}stholz}, {\em Logarithmic forms and group varieties}, 
{ J.~reine angew.~Math.} {\bf 442} (1993), 19-62. 
%
\bibitem{Berczes2} { A.~Berczes, B.~Brindza and L.~Hajdu}, {\em On the power values of polynomials}, 
{ Publ.~Math.~Debrecen} {\bf 53} (1998), 375-381.
%
\bibitem{Berczes} { A.~B\'{e}rczes and I.~Pink}, {\em On the Diophantine equation $x^2+q^{2k}=y^n$},  
{ Archive Math.~(Basel)}, {\bf 91} (2008), 505-517.
%
\bibitem{BH} { Y.~Bilu and G.~Hanrot}, {\em Solving Thue equations of high degree},
{ J.~Number Th.}, {\bf 60} (1996), 373-392.
%
\bibitem{PrimDiv} { Y.~Bilu, G.~Hanrot and P.M.~Voutier}, {\em Existence of primitive
divisors of Lucas and Lehmer numbers. With an appendix by M.Mignotte},
{ J.~reine angew.~Math.} {\bf 539} (2001), 75-122.
%
\bibitem{Bosma} { W.~Bosma, J.~Cannon and C.~Playoust}, {\em The Magma Algebra System I. The user language}, 
{ J.~Symbolic Comput.} {\bf 24} (3-4) (1997), 235-265.
%
\bibitem{Bugeaud2} { Y.~Bugeaud, M.~Mignotte and S.~Siksek}, {\em Classical and modular approaches to 
exponantial Diophantine equations II. The Lebesque-Nagell equation}, 
{ Compositio Math.} {\bf 142} (2006), 31-62.
%
\bibitem{magma-handbook} { J.J.~Cannon and W.~Bosma} (Eds.) 
{\em Handbook of Magma Functions, Edition 2.13} (2006), 4350 pages.
%
\bibitem{CDLPS} { I.N.~Cangul, M.~Demirci, F.~Luca, A.~Pint\'{e}r and G.~Soydan}, 
{\em On the Diophantine equation $x^{2}+2^{a}11^{b}=y^{n}$}, { Fibonacci Quart.}, to appear.
%
\bibitem{CDIPS} { I.N.~Cangul, M.~Demirci, I.~Inam, F.~Luca and G.~Soydan}, 
{\em On the Diophantine equation $x^{2}+2^{a}3^{b}11^{c}=y^{n}$}, submitted.

\bibitem{Carm} { R.~D.~Carmichael}, {\em On the numerical factors of the arithmetic forms $\al^n\pm\be^n$},
{ The Annals of Mathematics}, 2nd Ser. {\bf 15}, No 1/4 (1913-1914), 30-48.
%
\bibitem{Cohn2}  { J.H.E.~Cohn}, {\em The Diophantine equation $x^{2}+2^{k}=y^{n}$}, 
{ Arch.~Math(Basel)} {\bf 59} (1992), 341-344.
%
\bibitem{Cohn1} { J.H.E.~Cohn}, {\em The Diophantine equation $x^{2}+C=y^{n}$}, 
{ Acta Arith.} {\bf 65}, No 4, (1993), 367-381.
%
\bibitem{Luca5} { E.~Goins, F.~Luca and A.~Togbe}, {\em On the Diophantine equation 
$x^2+2^{\al }5^{\be}13^{\ga}=y^n$}, ANTS VIII Proceedings: A.~van der Poorten and A.~Stein (eds.), 
{ ANTS VIII, Lecture Notes in Computer Science} {\bf 5011} (2008), 430-442.
%
\bibitem{Gyory} { K.~Gy\H{o}ry, I.~Pink and A.~Pint\'{e}r}, {\em Power values of
polynomials and binomial Thue-Mahler equations}, { Publ.~Math.~Debrecen} {\bf 65} (2004), 341-362.
%
\bibitem{Landau} { E.~Landau and A.~Ostrowski}, {\em On the Diophantine equation $ay^2+by+c=dx^n$},  
{ Proc.~London Math.~Soc.} {\bf 19}, No 2 (1920), 276-280.
%
\bibitem{Le1} { M.H.~Le}, {\em On Cohn's conjecture concerning the Diophantine equation 
$x^{2}+2^{m}=y^{n}$},  { Arch.~Math(Basel)} {\bf 78} (2002), 26-35.
%
\bibitem{Lebesque} { V.A.~Lebesgue}, {\em Sur l'impossibilit\'{e} en nombres
entieres de l' \'{e}quation $x^{m}=y^{2}+1$}, { Nouvelles Ann.~Math.} {\bf 9}, No 1 (1850), 178-181.
%
\bibitem{Liqun1} { T.~Liqun}, {\em On the diophantine equation $X^2+3^m=Y^n$}, 
{ Integers: Electronic J.~Combinatorial Number Theory} {\bf 8} (2008), 1-7.
%
\bibitem{Liqun2} { T.~Liqun}, {\em On the diophantine equation $x^2+5^m=y^n$}, { Ramanujan J.} {\bf 19}
(2009), 325-338.
%
\bibitem{Luca3} { F.~Luca}, {\em On a diophantine equation}, 
{ Bull.~Aus.~Math.~Soc.} {\bf 61} (2000), 241-246.
%
\bibitem{Luca1} { F.~Luca}, {\em On the equation $x^2+2^a3^b=y^n$},  
{ Int.~J.~Math.~Math.~Sci.} {\bf 29}, No 4 (2002), 239-244.
%
\bibitem{LucaReport}{ F.~Luca}, {\em Effective Methods for Diophantine Equations}, Winter School on
Explicit Methods in Number Theory, January 26-30, 2009.
%
\bibitem{Luca4} { F.~Luca and A.~Togbe}, {\em On the Diophantine equation $x^2+7^{2k}=y^n$},  
{ Fibonacci Quart.} {\bf 54} No 4 (2007), 322-326.
%
\bibitem{Luca2} { F.~Luca and A.~Togbe}, {\em On the Diophantine Equation $x^2+2^a5^b=y^n$},  
{ Int.~J.~Number Th.}, {\bf 4} No 6 (2008), 973-979.
%
\bibitem{Mignotte} { M.~Mignotte and B.M.M.~de Weger}, {\em On the Diophantine equations $x^2+74=y^5$ 
and $x^2+86=y^5$},  { Glasgow Math.~J.} {\bf 38} (1996), 77-85.
%
\bibitem{PZGH} { A.~Peth\H{o}, H.G.~Zimmer, J.~Gebel and E.~Herrmann},
{\em Computing all $S$-integral points on elliptic curves}, 
{ Math.~Proc.~Camb.~Phil.~Soc.} {\bf 127} (1999), 383-402.
%
\bibitem{Pink} { I.~Pink}, {\em On the Diophantine equation $x^{2}+2^{\al}3^{\be}5^{\ga}7^{\de}=y^{n}$},
{ Publ.~Math.~Debrecen} {\bf 70}, No 1-2 (2007), 149-166.
%
\bibitem{TdW1} { N.~Tzanakis and B.M.M.~de Weger}, {\em On the practical solution of
the Thue equation}, { J.~Number Th.}, {\bf 31} (1989), 99-132.
%
\bibitem{TdW2}  { N.~Tzanakis and B.M.M.~de Weger},  {\em How to explicitly solve a Thue-Mahler equation}, 
{ Compositio Math.} {\bf 84} (1992), 223-288.
%
\bibitem{Yu} { K.~Yu}, {\em Linear forms in $p$-adic logarithms II}, 
{ Compositio Math.} {\bf 74} (1990), 15-113.
%
\end{thebibliography}
\end{document}